\documentclass{amsart}

\usepackage[french,english]{babel}
\usepackage{amsmath,amssymb,amsthm}
\usepackage{stmaryrd}
\usepackage[all]{xy}
\usepackage{enumerate}

\theoremstyle{break}
\newtheorem{thm}{Th\'eor\`eme}[section]
\newtheorem{prop}[thm]{Proposition}
\newtheorem{lem}[thm]{Lemme}
\newtheorem{cor}[thm]{Corollaire}
\theoremstyle{definition}
\newtheorem{dfn}[thm]{Definition}
\newtheorem{rmq}{Remarque}[section]

\newtheorem{exe}{Exemple}[section]

\newcommand\C{\mathbb{C}}

\newcommand\Z{\mathbb{Z}}
\newcommand\N{\mathbb{N}}
\newcommand\Q{\mathbb{Q}}

\newcommand\Kcyc{K_{\operatorname {cyc}}}

\newcommand\Kpi{K_{\pi}}
\newcommand\OK{\mathcal{O}_K}

\newcommand\A{\mathbf A}
\newcommand\Ap{\mathbf A^+}
\newcommand\AL{\mathbf A_L}

\newcommand\Atilde{\widetilde{\mathbf A}}
\newcommand\Atildep{\widetilde{\mathbf A}^+}
\newcommand\AK{\mathbf A_K}
\newcommand\AKp{\mathbf A_K^+}

\newcommand\B{\mathbf B}

\newcommand\Bp{\mathbf B^+}
\newcommand\BK{\mathbf B_K}
\newcommand\BKlog{\mathbf B_{K, \log}}
\newcommand\BKlogp{\mathbf B_{K, \log}^+}
\newcommand\BL{\mathbf B_L}
\newcommand\BLp{\mathbf B_L^+}
\newcommand\BLlog{\mathbf B_{L, \log}}
\newcommand\BLlogp{\mathbf B_{L, \log}^+}
\newcommand\Blogp{\mathbf B_{\log}^+}

\newcommand\Btildep{\widetilde{\mathbf B}^+}
\newcommand\BKp{\mathbf B_K^+}
\newcommand\Acris{\operatorname{A_{cris}}}
\newcommand\Bcris{\operatorname{B_{cris}}}
\newcommand\BdR{\operatorname B_{\dR}}
\newcommand\BrigKp{\mathbf{B}_{\rig,K}^+}
\newcommand\Bst{\operatorname{B_{st}}}
\newcommand\E{\bold E}
\newcommand\EL{\bold E_L}

\newcommand\Etilde{\widetilde{\mathbf E}}
\newcommand\Etildep{\widetilde{\mathbf E}^+}

\newcommand\D{\operatorname{D}}

\newcommand\Dp{\operatorname{D}^+}
\newcommand\DL{\operatorname{D}_L}
\newcommand\DLp{\operatorname{D}_L^+}
\newcommand\Dcris{\operatorname{D_{cris}}}
\newcommand\Dlogp{\operatorname{D_{log}^+}}
\newcommand\Dlog{\operatorname{D_{log}}}
\newcommand\DlogLp{\D_{L}^{+,\uni}}
\newcommand\Dst{\operatorname{D_{st}}}

\newcommand\cris{\operatorname{cris}}
\newcommand\dR{\operatorname{dR}}
\newcommand\hf{\operatorname{hf}}

\newcommand\naif{\operatorname{naif}}
\newcommand\rig{\operatorname{rig}}
\newcommand\st{\operatorname{st}}

\newcommand\Rep{\operatorname{Rep}}

\newcommand\phiGamma{\ensuremath{(\varphi, \Gamma)}}
\newcommand\phiG{\ensuremath{(\varphi, G)}}
\newcommand\phiN{\ensuremath{(\varphi, N)}}

\newcommand\Proof{\textrm{\em Preuve : }}
\newcommand\Eproof{\hfill{{\large $\square$}}}

\newcommand\ssi{si et seulement si }

\newcommand\id{\operatorname{id}}
\newcommand\rg{\operatorname{rg}}

\newcommand\Fil[1]{\operatorname{Fil}^{#1}}
\newcommand\Filhat[1]{\widehat{\operatorname{Fil}}{}^{#1}}
\newcommand\Frac{\operatorname{Frac }}
\newcommand\Gal{\operatorname{Gal }}
\newcommand\Hom{\operatorname{Hom}}

\newcommand\Sym{\operatorname{Sym}}
\newcommand\NN{\mathcal N}
\newcommand\Nlog{\mathcal N_{\log}}
\newcommand\uni{\operatorname{uni}}

\title[Repr\'esentations na\"ives et modules de Wach]{Sur les repr\'esentations absolument na\"ives et leur module de Wach}

\author{Floric Tavares Ribeiro}
\email{floric.tavares-ribeiro@unicaen.fr}
\address{Laboratoire de Math\'ematiques Nicolas Oresme \\
Universit\'e de Caen \\ 
CNRS UMR 6139 \\
Campus II\\
Boulevard Mar\'echal Juin \\
B.P. 5186 \\
14032 Caen Cedex\\
France}
\subjclass[]{11F80}
\keywords{repr\'esentations $p$-adiques,  ($\varphi,\Gamma$)-modules, repr\'esentations semi-stables, modules de Wach}

\date{}

\begin{document}

\maketitle

\selectlanguage{english}
\begin{abstract}
We construct a Wach module for the absolutely semi-stable representations  the filtered \phiN-module of which satisfies the Griffiths transversality, which happens in particular for ordinary representations. This construction extends the ones of Wach and Berger for the absolutely crystalline case. We use it in particular to describe semi-stable representations which factor through the Tate curve extension.
\end{abstract}

\selectlanguage{french}
\begin{abstract}
On construit un module de Wach pour les repr\'esentations absolument semi-stables dont le \phiN-module filtr\'e satisfait la transversalit\'e de
Griffiths, ce qui est le cas en particulier des repr\'esentations ordinaires. Ceci prolonge les constructions de Wach et Berger pour le cas absolument cristallin. On en d\'eduit notamment une description des repr\'esentations semi-stables se factorisant par l'extension de la courbe de Tate.
\end{abstract}

\tableofcontents

\section{Introduction}
\subsection{Notations}
Soit $p$ un nombre premier et $k$ un corps parfait de caract\'eristique $p$. On note $\OK$ l'anneau des vecteurs de Witt sur $k$ et $K$
son corps des fractions.
On fixe $\overline K$ une cl\^oture alg\'ebrique de $K$ et on note $G_K=\Gal(\overline K/K)$ le groupe de Galois absolu de $K$. On note $v_p$ la valuation de $\overline K$ normalis\'ee par $v_p(p)=1$ et $\C_p$ le compl\'et\'e de $\overline K$ pour $v_p$.
On fixe encore $\varepsilon=(\zeta_{p^n})_{n\in \N}$ une famille coh\'erente de racines primitives $p^n$-i\`emes de l'unit\'e :
$$ \zeta_1=1 \ ; \ \zeta_p \neq 1 \ ; \ \zeta_{p^{n+1}}^p=\zeta_{p^n} \forall n \in \N.$$
Le groupe de Galois $\Gamma = \Gal(\Kcyc/K)$ o\`u $\Kcyc$ est l'extension cyclotomique $\Kcyc = \cup_{n\in \N}K(\zeta_{p^n})$ s'identifie \`a $\Z_p^*$
par le caract\`ere cyclotomique $\chi$, de noyau $H_K = \Gal(\overline K/\Kcyc)$.

On appelle $\Z_p$-repr\'esentation (de $G_K$) tout $\Z_p$-module libre de type fini muni d'une action lin\'eaire et continue de $G_K$.
Une repr\'esentation $p$-adique (de $G_K$) est un $\Q_p$-espace vectoriel de dimension finie avec une action lin\'eaire et continue de $G_K$.
On obtient une repr\'esentation $p$-adique \`a partir d'une $\Z_p$-repr\'esentation en \'etendant les scalaires de $\Z_p$ \`a $\Q_p$, et
r\'eciproquement, toute repr\'esentation $p$-adique contient un $\Z_p$-r\'eseau stable par $G_K$ qui est donc une $\Z_p$-repr\'esentation,
que l'on appelle un r\'eseau de la repr\'esentation.

\subsection{Modules de Wach et repr\'esentations na\"ives}
Le module de Wach d'une repr\'esentation absolument cristalline \'etablit un pont entre les deux pendants
de la th\'eorie de Fontaine : le monde des \phiGamma-modules et celui de la th\'eorie de Hodge
$p$-adique. Wach \cite{Wach_pot_cris} a ainsi montr\'e qu'une repr\'esentation de hauteur finie \'etait cristalline
\`a la condition n\'ecessaire et suffisante de l'existence dans son \phiGamma-module d'un sous-module satisfaisant des conditions
techniques assez fortes.

Berger a pr\'ecis\'e les r\'esultats de Wach dans \cite{Berger_limitesrepcris} o\`u il d\'efinit la cat\'egorie des modules de Wach qu'il
prouve \'equivalente \`a celle des repr\'esentations cristallines. Pour cela, il ajoute une condition technique
que doit satisfaire le module de Wach $\NN$, ce qui le rend unique. Il fournit alors une filtration
sur ce module qui permet d'identifier le quotient $\NN/X\NN$ avec le module cristallin $\Dcris(V)$ comme $\varphi$-modules filtr\'es.
On renvoie au paragraphe \ref{ParagCasClassique} pour les d\'etails.

La th\'eorie de Hodge $p$-adique fait correspondre \`a une repr\'esentation semi-stable $V$ un $\phiN$-module filtr\'e $\Dst(V)$.
En g\'en\'eral, on n'attend pas que la monodromie $N$ respecte la filtration. On s'int\'eresse ici aux repr\'esentations semi-stables, dites na\"ives,
pour lesquelles la filtration v\'erifie la transversalit\'e de Griffiths,
c'est-\`a-dire $N(\Fil i) \subset \Fil {i-1}$ ;
c'est le cas en particulier des repr\'esentations ordinaires.

On associe de mani\`ere naturelle \`a une repr\'esentation na\"ive $V$ une repr\'esentation cristalline.
On en d\'eduit la construction d'un module de Wach $\NN(V)$ dans le \phiG-module $D_L(V)$ adapt\'e \`a l'extension de la courbe de Tate
$L = \cup_{n\in \N} K(\zeta_{p^n}, \sqrt[p^n]p)$ de groupe de Galois $G \simeq \overline{<\tau>} \rtimes \Gamma$, o\`u $\overline{<\tau>} \simeq \Z_p$.
On montre encore comment $\NN(V)/X\NN(V)$ s'identifie \`a $\Dst(V)$ comme \phiN-module filtr\'e, la filtration \'etant donn\'ee par la m\^eme
recette que celle de Berger, et la monodromie provenant de l'action de $\tau$.
Plus pr\'ecis\'ement, on obtient les \'enonc\'es suivants (Th\'eor\`emes \ref{Wachinitsst}, \ref{Wachsst} et \ref{Wachsstfiltre}
ci-dessous, on renvoie au texte pour les notations).

\begin{thm}
Si $V$ est une repr\'esentation $p$-adique de dimension $d$, alors $V$ est na\"ive si et seulement
s'il existe un sous-$\BKp$-module $\NN$ libre de rang $d$ de $\DLp(V)$, stable par $G$,
et un entier $h \in \Z$ tels que l'action de $\Gamma$ sur $(\NN/X\NN)(h)$ est triviale.
\end{thm}

\begin{thm}
Si $V$ est une repr\'esentation na\"ive positive de dimension $d$, alors il existe un unique
sous-$\BKp$-module $\NN(V)$ de $\DLp(V)$ satisfaisant les conditions suivantes :
\begin{enumerate}
\item $\NN(V)$ est un $\BKp$-module libre de rang $d$ qui engendre $\DL(V)$;
\item l'action de $G$ pr\'eserve $\NN(V)$, et est triviale sur $\NN(V)/X\NN(V)$ ;
\item il existe un entier $r\geq 0$ tel que $X^r\DlogLp(V) \subset \NN(V)$
\end{enumerate}
o\`u $\DlogLp(V)$ est le plus grand sous-$\BKp$-module de $\DLp(V)$ sur lequel $\tau$ est unipotent.
De plus, $\NN(V)$ est stable par $\varphi$.

Enfin, la somme $\frac 1 X \log \tau$ converge vers un op\'erateur $N$ nilpotent sur $\NN(V)$.
Si l'on munit de plus  $\NN(V)$ de la filtration
$$ \Fil i \NN(V) = \left\{x \in \NN(V), \varphi(x) \in \left(\frac{\varphi(X)}X\right)^i \NN(V)\right\},$$
alors il existe un isomorphisme
$\Dst(V) \simeq \NN(V)/X\NN(V)$ de $\phiN$-modules filtr\'es.
\end{thm}

Ceci nous permet de montrer que les repr\'esentations semi-stables se factorisant par $G$ sont na\"ives ainsi qu'une description de ces repr\'esentations (Proposition \ref{descriptionsstG}). Comme cons\'equence, on obtient une preuve dans le cas non ramifi\'e du th\'eor\`eme suivant d\^u \`a Kisin.

\begin{thm}
Le foncteur restriction de la cat\'egorie des repr\'esentations cristallines de $G_K$ vers la cat\'egorie des repr\'esentations de $\Gal(\Kpi/K)$  o\`u $\Kpi=\cup_{n\in \N} K(\sqrt[p^n]p)$ est pleinement fid\`ele.
\end{thm}
Remarquons que l'on peut tirer de ce travail les ingr\'edients d'une preuve courte du th\'eor\`eme de Kisin \cite{BeilinsonTavares}.

\bigskip

\section{Th\'eorie de Fontaine}
\subsection{\phiGamma-modules}
Commen\c cons par quelques rappels sur la th\'eorie des \phiGamma-modules. On renvoie le lecteur \`a l'article de Fontaine \cite{Fontaine_rep_p_corloc} pour les d\'etails, ainsi qu'\`a \cite{Tavares_AIF} pour la partie concernant les \phiGamma-modules adapt\'es \`a l'extension $L$.
\subsubsection{\phiGamma-modules usuels}
Soit $\Etilde$ la limite projective $\lim_{\leftarrow} \C_p$ o\`u les applications de transition sont l'\'el\'evation \`a la puissance $p$.
C'est un corps de caract\'eristique $p$, alg\'ebriquement clos et complet pour la valuation $v_\E((x_n)_n)=v_p(x_0)$.
Son anneau d'entiers $\Etildep$ s'identifie \`a $\lim_{\leftarrow} \mathcal O_{\C_p}$. Le syst\`eme $\varepsilon$ d\'efinit un \'el\'ement de $\Etildep$ encore not\'e $\varepsilon$ et on d\'efinit $\E_K = k((\varepsilon -1))$ et $\E$ la cl\^oture s\'eparable
de $\E_K$ dans $\Etilde$.\\
Le corps $\Etilde$ est parfait et l'on consid\`ere $\Atilde = W(\Etilde)$ l'anneau des vecteurs de Witt \`a coefficients dans
$\Etilde$ et $\Atildep$ le sous-anneau des vecteurs de Witt \`a coefficients dans $\Etildep$. Le corps $\E_K$ se rel\`eve en caract\'eristique $0$
en $\A_K$, le compl\'et\'e $p$-adique de $\OK[[X]][\frac 1 X]$ o\`u $X=[\varepsilon]-1$. C'est l'anneau des s\'eries de Laurent
$\sum a_nX^n$ avec pour tout $n\in \Z$,  $a_n \in \OK$ et $a_n$ tendant $p$-adiquement vers $0$ en $-\infty$.
On note $\A$ le compl\'et\'e $p$-adique de la cl\^oture non ramifi\'ee de $\A_K$ dans $\Atilde$. On note encore $\Ap = \Atildep \cap \A$ et
$\AKp = \Atildep \cap \A_K = \OK[[X]]$. Comme $\E^{H_K} = \E_K$, on a aussi $\A^{H_K} = \A_K$ et de m\^eme avec les versions avec un plus.
L'action r\'esiduelle de $\Gamma$ et l'action induite par le Frobenius $\varphi$ de $\Atilde$ sur $\A_K$, qui prolonge le Frobenius naturel de $K$,
se d\'ecrivent alors ainsi :
$$\varphi(X) = (1+X)^p-1 \\ \ ; \ \ \ \sigma(X) = (1+X)^{\chi(\sigma)}-1, \ \sigma \in \Gamma.$$
On d\'efinit \`a partir de ces anneaux les corps $\B=\A\otimes_{\Z_p}\Q_p$, $\B_K=\A_K\otimes_{\Z_p}\Q_p$, et de m\^eme pour
les anneaux $\Bp$, $\BKp$ et les versions avec tilde.

Un $\phiGamma$-module $D$ sur $\A_K$ (resp. $\BK$) est un module libre de type fini, muni d'actions continues semi-lin\'eaires commutant de $\varphi$ et $\Gamma$.
Il est \'etale si $\varphi(D)$ engendre $D$ sur $\A_K$ (resp. s'il poss\`ede un $\AK$-r\'eseau \'etale).

Rappelons le th\'eor\`eme fondamental des $\phiGamma$-modules.
\begin{thm}[Fontaine \cite{Fontaine_rep_p_corloc}]
Le foncteur  $ \D : T \mapsto \D(T) = (\A\otimes T)^{H_K}$ (respectivement $ \D : V \mapsto \D(V) = (V \otimes \B)^{H_K}$)
d\'efinit une \'equivalence de cat\'egories tannakiennes entre la cat\'egorie
des $\Z_p$-repr\'esentations (resp. repr\'esentations $p$-adiques) de $G_K$ et la cat\'egorie des $\phiGamma$-modules \'etales sur $\A_K$
(resp. $\BK$).

\end{thm}
Une $\Z_p$-repr\'esentation $T$ est dite de hauteur finie d\`es lors que $\Dp(T) = (\Ap\otimes T)^{H_K}$ contient une base de $\D(T)$ sur $\A_K$.
Une repr\'esentation $p$-adique $V$ est de hauteur finie si elle contient un r\'eseau de hauteur finie.
De mani\`ere \'equivalente, un $\phiGamma$-module est de hauteur finie s'il contient une base dans laquelle la matrice de $\varphi$ est \`a
coefficients dans $\AKp$.

\subsubsection{($\varphi$, $G$)-modules}
On fixe $\rho = (\pi_n)_{n\in \N}$ une famille coh\'erente de racines $p^n$-i\`eme de $p$ : $\pi_0=p$ et $\pi_{n+1}^p=\pi_n$ pour tout $n$.
Alors $G_K$ agit sur $\rho$ par $g(\rho) = \rho \varepsilon^{\xi(g)}$ ce qui d\'efinit l'application $\xi : G_K \rightarrow \Z_p$, ou plus
pr\'ecis\'ement le $1$-cocycle $\xi\otimes\varepsilon$ \`a valeurs dans $\Z_p(1)$. \\
On d\'efinit la courbe de Tate $T_2$ comme le $\Z_p$-module $\Z_p \oplus \Z_pe$ de base $1, e$ muni de l'action
$$g(1)= 1 \ \ ; \ \ \ g(e) = \frac 1 {\chi(g)} (e + \xi(g)), \ g \in G_K.$$
Ainsi, $T_2(1)$, comme tout module de Tate d'une courbe elliptique \`a mauvaise r\'e\-duc\-tion multiplicative d\'eploy\'ee appara\^it comme une extension
$$0 \longrightarrow \Z_p(1) \longrightarrow T_2(1) \longrightarrow \Z_p \longrightarrow 0.$$
L'action de $G_K$ sur $T_2$ se factorise par $G=\Gal(L/K)$ o\`u $L = \cup_{n\in\N}K(\zeta_{p^n}, \pi_n)$ est le compositum de $\Kcyc$ et de l'extension
kummerienne $\Kpi = \cup_{n\in \N} K(\pi_n)$ (c'est aussi la cl\^oture galoisienne de $\Kpi$). Le groupe $G$ s'identifie \`a $\Z_p \rtimes \Gamma$
si l'on identifie $\Gamma$ \`a $\Gal(L/\Kpi)$ par $\chi$ et $\Z_p$ \`a $\Gal(L/\Kcyc)$ par $\xi$. Ceci fournit un g\'en\'erateur
topologique $\tau \in \Gal(L/\Kcyc)$ tel que $\Gal(L/\Kcyc)$ s'identifie \`a $\tau^{\Z_p}$. On a donc $\tau (\rho) = \rho \varepsilon$,
\emph{i.e.} $\xi(\tau)=1$.

$$\xymatrix{
& L \ar@{-}[dd]^{G} \ar@{-}[dl]_{\tau^{\Z_p}} \ar@{-}[dr]^{\Gamma} & \\
\Kcyc \ar@{-}[dr]_{\Gamma} & & \Kpi \ar@{-}[dl]\\
& K &}$$

 Notons que le choix de $\tau$ ne d\'epend pas de celui de $\rho$ mais bien du choix de $\varepsilon$. Changer $\varepsilon$
en $\varepsilon^\alpha$ revient \`a changer $\tau$ en $\tau^\alpha$. La relation de commutation entre $\Gamma$ et $\tau$ est donn\'ee par
$\sigma \tau \sigma^{-1} = \tau^{\chi(\sigma)}.$

\bigskip

De m\^eme que dans le cas cyclotomique, on peut d\'efinir le foncteur $$\DL(T) = (T \otimes \A)^{G_L}$$ qui \`a une $\Z_p$-repr\'esentation $T$ associe
un $\phiG$-module \'etale sur $\AL = \A^{G_L}$, \emph{i.e.} un $\AL$-module libre de type fini muni d'actions semi-lin\'eaires commutant de $\varphi$ et $G$ et $\varphi(\DL(T))$ engendrant $\DL(T)$ comme $\AL$-module. On d\'efinit de m\^eme les versions rationnelles de $\DL$ et de $\phiG$-module sur $\BL$.
On a encore (cf. \cite{Tavares_AIF}) :
\begin{thm}
Le foncteur $\DL$ induit une \'equivalence de cat\'egories tannakiennes entre la cat\'egorie
des $\Z_p$-repr\'esentations (respectivement re\-pr\'e\-sen\-ta\-tions $p$-adiques) de $G_K$ et la cat\'egorie des $\phiG$-modules \'etales sur $\AL$
(resp. $\BL$).
De plus $\DL \simeq \AL \otimes_{\A_K} \D$.
\end{thm}

\subsection{Th\'eorie de Hodge $p$-adique}
\subsubsection{Repr\'esentations semi-stables}
Ce paragraphe est consacr\'e \`a quelques rappels sur les p\'eriodes $p$-adiques des repr\'esentations semi-stables. Pour les constructions et d\'emonstrations originelles, on renvoie \`a \cite{Fontaine_cor_perio}.

Un \'el\'ement $x$ de $\Atildep$ s'\'ecrit de mani\`ere unique 
$$x = \sum_{n\in \N} p^n[x_n] \textrm{ avec } x_n = \left(x_n^{(m)}\right)_{m\in \N} \in \Etildep.$$
L'application $\theta : x \mapsto \sum_{n\in \N} p^nx_n^{(0)}$ est un morphisme surjectif de $\Atildep$ sur $\mathcal O_{\C_p}$ de noyau l'id\'eal principal engendr\'e par $\omega = X/\varphi^{-1}(X)$.
L'anneau $\BdR^+$ est le compl\'et\'e de $\Btildep$ pour la topologie $\omega$-adique. La s\'erie $\log (1+X)$ converge
dans $\BdR^+$. Sa limite est not\'ee $t$, et $\BdR = \BdR^+[1/t] = \Frac \BdR^+$.\\
On appelle $\Acris$ le compl\'et\'e $p$-adique de l'enveloppe \`a puissances divis\'ees de $\Atildep$ par rapport \`a $\ker \theta = (\omega)$.
Il s'identifie au sous-anneau de $\BdR^+$ form\'e des \'el\'ements de la forme $\sum_{n \in \N} a_n \frac {\omega^n}{n!}$ avec $a_n \in \Atildep$ tendant
$p$-adiquement vers $0$. La s\'erie d\'efinissant $t$ converge dans $\Acris$. On d\'efinit $\Bcris = \Acris[1/t]$.\\
Comme $\rho$ est un syst\`eme coh\'erent de racines $p^n$-i\`emes de $p$, il d\'efinit un \'el\'ement de $\Etildep$ et si l'on note
$Y = [\rho]$, la s\'erie $\log Y = \log \frac Y p + \log p$ converge, d\`es que l'on a choisi la valeur de $\log p$
(\emph{e.g.} $\log p = 0$), dans $\BdR^+$,  mais pas dans $\Bcris$. On note $\Bst = \Bcris[\log Y]$.\\
L'action de $G_K$ s'\'etend par continuit\'e \`a $\BdR$ ainsi qu'\`a $\Bcris$ et $\Bst$. De plus, $\BdR$ est filtr\'e par
$\Fil i \BdR = t^i \BdR^+$ ; on munit $\Bcris$ et $\Bst$ de la filtration induite. Ces deux derniers anneaux sont de plus
munis d'une action de $\varphi$ \'etendue par continuit\'e et la $\Bcris$-d\'erivation $N = -\frac d {d \log Y}$ fournit un op\'erateur localement nilpotent sur $\Bst$ tel que $\Bcris= (\Bst)^{N=0}$.

\bigskip

\`A toute repr\'esentation $p$-adique $V$ on associe le $K$-espace vectoriel $$\Dst(V) = (V \otimes \Bst)^{G_K}$$ de dimension inf\'erieure ou \'egale
\`a $\dim V$. Quand il y a \'egalit\'e, on dit que $V$ est semi-stable. Les structures de $\Bst$ induisent sur $\Dst(V)$ celle d'un  $\phiN$-module filtr\'e, \emph{i.e.} un $K$-espace vectoriel muni d'une filtration $\Fil i $ par des sous-$K$-espaces vectoriels, d\'ecroissante, exhaustive et s\'epar\'ee, ainsi que
d'une application semi-lin\'eaire injective $\varphi$ et d'une application $K$-lin\'eaire nilpotente $N$ v\'erifiant
$N\varphi = p\varphi N$. Soit $t_H(D)$ l'unique saut de la filtration de $\det D$ et $t_N(D)$ la valuation $p$-adique de $\varphi$ sur $\det D$.
On dit qu'un $\phiN$-module filtr\'e $D$ est admissible si $t_H(D) = t_N(D)$ et $t_H(D') \leq t_N(D')$ pour tout sous-objet $D'$ de $D$.
\begin{thm}[Fontaine, Colmez-Fontaine]
Le foncteur $\Dst$ d\'efinit une \'equi\-va\-len\-ce de cat\'egories tannakiennes entre la cat\'egorie $\Rep_{\st} G_K$ des re\-pr\'e\-sen\-ta\-tions
semi-stables de $G_K$ et la cat\'egorie des $\phiN$-modules filtr\'es admissibles.
\end{thm}
La partie difficile de ce th\'eor\`eme est la preuve que tout module admissible provient d'une repr\'esentation semi-stable. Des preuves en ont \'et\'e
donn\'ees par Colmez-Fontaine \cite{ColmezFontaine_ConstructionSst}, puis Colmez \cite{Colmez_Banachfinie}, Fontaine \cite{Fontaine_ArithmetiqueRepGalois},
Berger \cite{Berger_EquadiffModulesfiltres}, et Kisin \cite{Kisin_crys}.

Une repr\'esentation semi-stable est dite cristalline si elle est semi-stable et $N=0$ sur son $\phiN$-module filtr\'e, ou de mani\`ere \'equivalente
si $\Dcris(V) = \Dst(V)^{N=0} = (V \otimes \Bcris)^{G_K}$ est de dimension maximale $\dim V$.
De plus, une repr\'esentation $V$ semi-stable est dite positive d\`es lors que $\Fil 0 \Dst(V) = \Dst(V)$. Pour des raisons de confort, on supposera souvent les repr\'esentations positives ;
cette condition n'est pas restrictive car toute repr\'esentation semi-stable peut \^etre rendue positive en la tordant par une puissance (n\'egative)
suffisante du caract\`ere cyclotomique.

\subsubsection{La courbe de Tate et $\Bst$}
La repr\'esentation $V_2 = \Q_p\otimes_{\Z_p}T_2$ associ\'ee \`a la courbe de Tate est semi-stable, non cristalline :
$\Dst(V_2)$ admet pour base $ \epsilon_1 = 1\otimes 1$, $\epsilon_2 =  1\otimes \log Y -  e\otimes t$. On a donc
$$\varphi(\epsilon_1) = \epsilon_1,\ \varphi(\epsilon_2) = p\epsilon_2 \ \textrm{  et } \
N(\epsilon_2)=\epsilon_1, \ N(\epsilon_1)=0.$$
En fait, on remarque que $\Z_p + \frac {\log Y} t \Z_p$ s'identifie \`a $T_2 = \Z_p + e \Z_p$ comme $\Z_p$-repr\'esentation de $G_K$.\\
De m\^eme, les repr\'esentations $V_d = \Q_p\otimes_{Z_p}T_d$ associ\'ees \`a $T_d = \Sym^{d-1}T_2$ pour $d\geq 2$ sont semi-stables de dimension $d$ ;
on peut identifier $T_d$ au groupe additif $\Z_p\left[\frac {\log Y} t\right]_{<d}$ des polyn\^omes en $\frac {\log Y} t$ de degr\'e plus petit que $d-1$ muni de l'action naturelle de $G_K$. Ces repr\'esentations forment ainsi un syst\`eme inductif, et \`a la limite :
\begin{prop}\label{descriptionBst}
Si $\displaystyle{V_\infty = \lim_{\rightarrow} \Sym^d V_2}$, alors  $V_\infty \simeq \Q_p\left[e\right] $ et l'application $\Bcris$- lin\'eaire
$\Bcris\otimes_{\Q_p}V_\infty \rightarrow \Bst$ qui associe $\frac {\log Y} t$ \`a $e$ est un isomorphisme compatible
avec l'action de $G_K$ et de $\varphi$. L'op\'erateur $N$ co\"incide via cette identification \`a $-\frac 1 t \frac d {de}$.
\end{prop}
\Proof Il suffit de remarquer que $ \varphi\left(\frac {\log Y} t\right) = \frac {\log Y} t$. \Eproof

On munit $\Bst$ de la filtration $\Filhat {}$ induite par celle de $\Bcris$ sur $\Bcris\otimes_{\Q_p}V_\infty$ :
$$\Filhat i \Bst = \left(\Fil i \Bcris\right)\otimes_{\Q_p}V_\infty.$$
Cette filtration ne co\"incide pas avec la filtration $\Fil {}$ de $\Bst$ comme le montre la proposition suivante.

\begin{prop}\label{filhat}
Pour tout $i \in \Z$, $$ \Filhat i \Bst = \bigcap_{k \in \N} N^{-k}\left(\Fil {i-k} \Bst\right).$$
En particulier $N\left(\Filhat i \Bst\right) \subset \Filhat {i-1} \Bst$ pour tout $i \in \Z$, c'est-\`a-dire que $\Filhat {}$ respecte la
transversalit\'e de Griffiths.
\end{prop}
\Proof
Pour tout $n\ge 0$ on pose $E_n = \ker N^n$. Alors  pour tous $0 \le k\le n$, $N^k(E_n)= E_{n-k}$ et l'on sait que $\Bst = \cup_{n \in \N} E_n$ et $E_0 = \Bcris$. On proc\`ede par r\'ecurrence sur $n$ et l'on note
$\Filhat i E_n = \Filhat i \Bst \cap E_n$. On veut donc montrer que
$$ \forall i, n \in \N, \ \Filhat i E_n = \bigcap_{k =0}^n N^{-k}\left(\Fil {i-k} E_{n-k}\right),$$
c'est-\`a-dire que si $x\in E_n$, alors
$$x\in \Filhat i E_n \Longleftrightarrow  \forall 0\leq k \leq n, \ N^k(x) \in \Fil {i-k} E_{n-k}.$$
Ceci est imm\'ediat pour $E_0 = \Bcris$. Et si $x\in E_{n+1}$,
on l'\'ecrit $$x = \sum_{m=0}^n a_m\left(\frac {\log Y} t\right)^m.$$
Ainsi $N(x) = \sum_{m=1}^n \frac {ma_m} t \left(\frac {\log Y} t\right)^{m-1}$.
On a donc
\begin{eqnarray*}
x\in \Filhat i E_n  & \Longleftrightarrow & \forall 0 \leq m \leq n, \ a_m \in \Fil i \Bcris\\
& \Longleftrightarrow & a_0 \in \Fil i \Bcris \textrm{ et } N(x) \in \Filhat {i-1} \Bst\\
& \Longleftrightarrow & x \in \Fil i \Bst \textrm{ et } N(x) \in \Filhat {i-1} \Bst
\end{eqnarray*}
ce qui permet de conclure par r\'ecurrence. La derni\`ere \'equivalence provient de ce que $\frac {\log Y} t \in \Fil 0 \Bst$.
\Eproof

Breuil construit dans \cite[\S 2]{Breuil_Griffiths}  un anneau $\widehat {\operatorname B}_{\st}$, compl\'et\'e de $\Bst$, muni d'une filtration
naturelle respectant la transversalit\'e de Griffiths. La filtration qu'elle induit sur $\Bst$ co\"incide avec celle que l'on a d\'efinie ici.

Suivant Breuil, on dira qu'un $\phiN$-module filtr\'e est na\"if s'il v\'erifie la transversalit\'e de Griffiths et, de m\^eme, qu'une repr\'esentation semi-stable est na\"ive si son $\phiN$-module est na\"if.
C'est le cas en particulier des repr\'e\-sen\-ta\-tions ordinaires \cite[Propositions 5.2. et 5.4.]{Mokrane_Ordinarite},
par exemple le module de Tate d'une courbe elliptique \`a mauvaise r\'eduction multiplicative d\'eploy\'ee.

\begin{dfn}
Si $\Fil {}$ est la filtration d'un $\phiN$-module filtr\'e on d\'efinit la filtration $\Filhat{}$ par
$$ \Filhat i= \bigcap_{k \in \N} N^{-k}\left(\Fil {i-k}\right),$$
\end{dfn}

Le lemme suivant est imm\'ediat.

\begin{lem}
Si $\Fil {}$ est la filtration d'un $\phiN$-module filtr\'e alors $\Filhat{}$ satisfait la transversalit\'e de Griffiths ; de plus 
$\Fil{}$  satisfait la transversalit\'e de Griffiths si et seulement si $\Fil{} = \Filhat{}$.
\end{lem}

Un \phiN-module filtr\'e est donc na\"if \ssi sa filtration $\Fil{}$ co\"incide avec $\Filhat{}$. On a \'egalement la caract\'erisation suivante.

\begin{lem}\label{filtrationsetadmissibilite}
Si $D$ est un $\phiN$-module filtr\'e admissible tel que le $\varphi$-module $D$ muni de la filtration $\Filhat{}$ est admissible, alors $D$ est na\"if.
\end{lem}
\Proof
Notons $D'$ le $\varphi$-module $D$ muni de la filtration $\Filhat{}$.
On fixe une base $(b_1, \dots, b_d)$ de $D$ adapt\'ee \`a la filtration $\Filhat{}$, c'est-\`a-dire qu'il existe une suite croissante d'entiers $(k_i)$
tels que pour tout $i$,
$\Filhat i D = \bigoplus_{k \geq k_i} Kb_k$. Pour tout $1 \leq k \leq d$,
il existe des entiers $i_k$ et $j_k$ tels que $b_k \in \Filhat {i_k} D \setminus \Filhat {i_k+1} D$ et
$b_k \in \Fil {j_k} D \setminus \Fil {j_k+1} D$ ; on a alors
$t_H(D') = \sum_{k=1}^d i_k$ et $t_H(D) \geq \sum_{k=1}^d j_k$. De $\Filhat i \subset \Fil i$ on  d\'eduit que
pour tout $k$, $i_k \leq j_k$, et donc, puisque $D$ et $D'$ sont admissibles,
$$t_N(D') = t_H(D') = \sum_{k=1}^d i_k \leq \sum_{k=1}^d j_k \leq t_H(D) = t_N(D)= t_N(D').$$
Ainsi, toutes les in\'egalit\'es sont des \'egalit\'es et $i_k=j_k$ pour tout $k$.

D\'eduisons-en maintenant que les filtrations co\"incident. On a l'inclusion $\Filhat i \subset \Fil i$ ; pour montrer la r\'eciproque, choisissons $x = \sum_{k =1}^d \alpha_k b_k \in \Fil i D$.
Pour tout $k$ tel que $\alpha_k \neq 0$ la famille $(b_1, \dots, b_{k-1}, x, b_{k+1}, \dots, b_d)$ est encore une base de $D$ et
$0 \neq b_1 \wedge \dots \wedge b_{k-1} \wedge x \wedge \dots \wedge b_d \in \Fil h (\det D)$ avec $h =i+\sum_{r\neq k} j_r$.
On en d\'eduit que $j_k \geq i$, et ceci \'etant vrai pour tout $k$ tel que $\alpha_k \neq 0$, $x = \sum_{k \geq k_i} \alpha_k b_k \in \Filhat {i}D$, comme attendu.
\Eproof

\begin{rmq}
On n'a en fait utilis\'e de la transversalit\'e de Griffiths que le fait que $\Filhat i \subset \Fil i$ et de l'admissibilit\'e de $D'$ que l'\'egalit\'e $t_N(D') = t_H(D')$. 

Nous verrons plus bas une r\'eciproque \`a ce lemme (Proposition \ref{cristallineassociee}).
\end{rmq}

\section{Modules de Wach}
\subsection{Le cas absolument cristallin}\label{ParagCasClassique}
Wach \cite{Wach_pot_cris} a donn\'e la ca\-rac\-t\'e\-ri\-sa\-tion suivante des repr\'esentations absolument cristallines de hauteur finie.
\begin{thm}[Wach]\label{Wach}
Si $V$ est une repr\'esentation $p$-adique de hauteur finie et de dimension $d$, alors $V$ est cristalline si et seulement
s'il existe un sous-$\BKp$-module libre $\NN$ de rang $d$ de $\Dp(V)$ stable par $\Gamma$ et un entier $h \in \Z$ tels que l'action de $\Gamma$
sur $(\NN/X\NN)(h)$ est triviale.
\end{thm}
Par ailleurs, Colmez \cite{ColmezCristallinesHauteurFinie} a montr\'e que toute repr\'esentation absolument cristalline est de hauteur finie si bien que l'on peut supprimer cette hypoth\`ese.
De plus, Berger \cite[Proposition II.1.1.]{Berger_limitesrepcris} a pr\'ecis\'e ce r\'esultat en ajoutant une condition qui rend le module de Wach unique.

\begin{thm}[Berger]\label{Berger1}
Si $T$ est un r\'eseau d'une repr\'esentation $V$ cristalline positive, alors il existe un
unique sous-$\AKp$-module $\NN(T)$ de $\Dp(T)$ satisfaisant les conditions suivantes :
\begin{enumerate}
\item $\NN(T)$ est un $\AKp$-module libre de rang $d = \dim_{\Q_p}V$ ;
\item l'action de $\Gamma$ pr\'eserve $\NN(T)$ et est triviale sur $\NN(T)/X\NN(T)$ ;
\item il existe un entier $r\geq 0$ tel que $X^r\Dp(T) \subset \NN(T)$.
\end{enumerate}
De plus, $\NN(T)$ est stable par $\varphi$. Enfin, si l'on pose $\NN(V) = \BKp\otimes_{\AKp}\NN(T)$, alors
$\NN(V)$ est l'unique sous-$\BKp$-module de $\Dp(V)$ qui v\'erifie l'analogue des conditions ci-dessus.
\end{thm}
Soit $\BrigKp$ l'anneau des s\'eries $\sum_{k=0}^{\infty}a_kX^k$ avec $a_k \in K$ qui convergent sur le disque unit\'e ouvert.
Il contient naturellement $\BKp$ comme sous-anneau. Berger a d\'ecrit \cite{Berger_equadiff} comment retrouver les modules $\Dst(V)$ et $\Dcris(V)$
\`a partir de $\phiGamma$-modules sur des anneaux de p\'eriodes surconvergentes proches de $\BrigKp$. En particulier dans le cas
d'une repr\'esentation $V$ positive cristalline il montre \cite[Proposition II.2.1.]{Berger_limitesrepcris} que $\Dcris(V) = (\BrigKp \otimes_{\BKp}\NN(V))^\Gamma$. Il en
d\'eduit la recette suivante \cite[Th\'eor\`eme III.4.4.]{Berger_limitesrepcris} qui permet de reconstruire $\Dcris(V)$ \`a partir de $\NN(V)$.
\begin{thm}[Berger]\label{Berger2}
Si $V$ est une repr\'esentation cristalline positive et si l'on munit $\NN(V)$ de la filtration
$$ \Fil i \NN(V) = \left\{x \in \NN(V), \varphi(x) \in \left(\frac{\varphi(X)}X \right)^i\NN(V)\right\},$$
alors l'inclusion $\Dcris(V) \subset \BrigKp \otimes_{\BKp}\NN(V)$ induit un isomorphisme
$$ \Dcris(V) \simeq \NN(V)/X\NN(V)$$ de $\varphi$-modules filtr\'es.
\end{thm}

\subsection{R\'eduction du cas na\"if au cas cristallin}
On construit un analogue du module de Wach en se ramenant au cas cristallin gr\^ace \`a la proposition suivante.
\begin{prop}\label{cristallineassociee}
Si $V$ est une repr\'esentation $p$-adique de dimension $d$, si l'on note $V[n] = V \otimes_{\Q_p} V_n$ pour $n \ge 2$, et si l'on munit
$\Dst(V)$ de la filtration $\Filhat{}$, alors
\begin{enumerate}[(i)]
\item l'application $e \mapsto \frac {\log Y} t$ induit pour tout $n\ge d$ un isomorphisme de $\varphi$-modules filtr\'es $\Dcris(V[n]) \simeq \Dst(V)$ ;
\item $V$ est na\"ive si et seulement si la plus grande sous-repr\'esentation cristalline $\widetilde V$ de $V[d]$ est de dimension $d$ ;
\item dans ce cas, $\Dcris(\widetilde V) = \Dcris(V[d]) $ est stable par l'op\'erateur $-\frac 1 t \frac d {de}$ qui induit la monodromie $N$ sur $\Dst(V)$.
\end{enumerate}
\end{prop}
\Proof
La Proposition \ref{descriptionBst} donne $(i)$.

On en d\'eduit encore que si la plus grande sous-repr\'esentation cristalline $\widetilde V$ de $V[d]$ est de dimension $d$, alors $\Dcris(V[d]) \simeq \Dst(V)$ est de dimension au moins $d$, donc $V$ est semi-stable
et $\Dcris(\widetilde V) \simeq \Dst(V)$ si $\Dst(V)$ est muni de la filtration $\Filhat{}$.
Le fait que $V$ est na\"ive est alors une cons\'equence du Lemme \ref{filtrationsetadmissibilite}.

Supposons $V$ na\"ive ; alors $V$ est semi-stable, et $V[d]$ l'est \'egalement.
Il s'agit de montrer
que le $\varphi$-module filtr\'e $\Dcris(V[d])$ de dimension $d$ est admissible.
Or $\Dcris(V[d])$ est un sous-objet du $\phiN$-module filtr\'e admissible $\Dst(V[d])$ et donc, tout sous-objet $D$ de
$\Dcris(V[d])$ est un sous-objet de $\Dst(V[d])$ et v\'erifie $t_H(D)\leq t_N(D)$.
De plus,
$$\det \Dst(V) \simeq \det \Dcris(V[d]))$$ comme $\varphi$-modules filtr\'es si bien que
$\Dcris(V[d])$ est admissible. Enfin, la description de l'op\'erateur $N$ sur $\Dcris(V[d])$ provient encore de  la Proposition \ref{descriptionBst}.
\Eproof

D\'efinissons $$\DLp(V)= (\Bp\otimes V)^{G}.$$
Nous allons construire un module interm\'ediaire entre $\Dp(V)$ et $\DLp(V)$ contenant le module de Wach de $\widetilde V$.
On d\'efinit $\Blogp = \Bp\otimes_{\Q_p}V_\infty = \Bp[e]$ muni de sa structure de $\Bp$-alg\`ebre, du Frobenius $\varphi = \varphi\otimes \id$ et de l'action naturelle de $G_K$. Rappelons que l'on a d\'efini l'application $\xi : G_K \rightarrow \Z_p$ par $g(\rho) = \rho \varepsilon^{\xi(g)}$.
\begin{lem}\label{elementb}
Il existe $b \in \BL$ unique modulo $\BK$ tel que pour tout $h \in H_K$, $(h-1)b = \xi(h)$ ; de plus $b$ est transcendant sur le compositum $\BK.\Frac\Bp$ et l'on a 
\begin{eqnarray*}
\left(\Bp[b]\right)^{H_K} = \left(\BLp[b]\right)^{\tau=1}=\BKp,\\
\left(\BK.\Bp[b]\right)^{H_K} = \left(\BK.\BLp[b]\right)^{\tau=1}=\BK.
\end{eqnarray*}
\end{lem}
\Proof Calculons une base du $\phiGamma$-module $\D(V_2)$. Le vecteur $1\otimes 1$ est stable par $H_K$.
Un vecteur non colin\'eaire s'\'ecrit $a\otimes e-b\otimes 1$ avec $a \neq 0$ si bien qu'en divisant on peut supposer
$a=1$. Ainsi, $b$ est bien uniquement d\'etermin\'e modulo $\BK$, et il satisfait la relation $(h-1)b = \xi(h)$ pour tout $h \in H_K$.

Remarquons d'abord que $b$ ne peut pas appartenir \`a $\Frac \Bp$ car sinon $V_2$ serait de hauteur finie et de de Rham, donc potentiellement cristalline d'apr\`es le th\'eor\`eme de Wach \cite{Wach_pot_cris}.

Supposons que $b \in \BK.\Frac\Bp$ et \'ecrivons $b= a_0 + \sum_{k=1}^n b_ka_k$ o\`u pour $0 \le k \le n$,  $a_k \in \Frac \Bp$ et pour $1 \le k \le n$, $b_k \in \BK$, et o\`u l'on suppose que $1, b_1, \dots, b_n$ sont lin\'eairement ind\'ependants sur $\Frac \Bp$.
Alors pour tout $h\in H_K$,
$$(h-1)b = (h-1)a_0 +  \sum_{k=1}^n b_k.(h-1)(a_k) = \xi(h) \in\Z_p$$
si bien que pour $1 \le k \le n$, $(h-1)a_k=0$, c'est-\`a-dire $a_k \in\BK$. 
L'\'el\'ement $b$ \'etant d\'efini modulo $\BK$, on peut ainsi choisir $b=a_0\in \Frac \Bp$, mais on a d\'ej\`a remarqu\'e que ceci \'etait exclu.

Supposons que $b$ (dont on fixe encore un repr\'esentant de la classe modulo $\BK$) soit alg\'ebrique sur le compositum $\BK.\Frac\Bp$ ; soit $P(X) =  X^n + \sum_{k=0}^{n-1} a_kX^k \in (\BK.\Frac\Bp)[X]$ le polyn\^ome minimal unitaire de $b$, de degr\'e $n$. Pour tout $h \in H_K$, on a donc :
$$ h(P(b)) =  (b+\xi(h))^n + \sum_{k=0}^{n-1} h(a_k)(b+\xi(h))^k =0.$$
On obtient ainsi pour tout $h \in H_K$ un polyn\^ome unitaire, annulateur de $b$ de degr\'e $n$ qui co\"incide par cons\'equent avec $P$, si bien que l'on a pour tout $h\in H_K$, 
$$h(a_{n-1}) +n\xi(h) = a_{n-1}.$$
Ainsi, $-a_{n-1}/n$ satisfait la propri\'et\'e d\'efinissant $b$, or on a montr\'e que $b$ n'appartenait pas \`a $\BK.\Frac\Bp$, d'o\`u il ressort que $b$ est transcendant sur ce corps.

Le m\^eme raisonnement montre que si $x=\sum_{k=0}^n a_kb^k \in \left(\BLp[b]\right)^{\tau=1}$ (respectivement $x \in (\BK.\BLp[b])^{\tau=1}$) avec $n\geq 1$ et $a_n \neq 0$, alors
$$ \tau(x) = \sum_{k=0}^n \tau(a_k)(b+1)^k = \sum_{k=0}^n a_kb^k $$
d'o\`u l'on d\'eduit que $a_n \in \BKp$ (respectivement $a_n \in \BK$) et $\tau(a_{n-1} + na_n) = a_{n-1}$, c'est-\`a-dire que 
$b+\frac {a_{n-1}}{na_n}\in \BK$ ce qui n'est pas possible. On en d\'eduit les \'egalit\'es.\Eproof
\begin{prop}\label{injectionlogb}
L'application $\BKp$-lin\'eaire $\iota : \Blogp \rightarrow \B$ (respectivement l'application $\BK$-lin\'eaire $\iota : \BK.\Blogp \rightarrow \B$) envoyant $e$ sur $b$ est injective et compatible \`a l'action de $H_K$.
\end{prop}
\Proof  C'est une cons\'equence du lemme pr\'ec\'edent, et du fait que l'action de $H_K$ sur $e$ est bien donn\'ee par $(h-1)e = \xi(h) = (h-1)b$. \Eproof

En combinant cette proposition avec le lemme pr\'ec\'edent, on d\'eduit encore l'\'egalit\'e $$\left(\Blogp\right)^{H_K}=\BKp.$$

\`A une repr\'esentation $p$-adique $V$ on associe le $\BKp$-module $$\Dlogp(V) = (\Blogp\otimes_{\Q_p}V)^{H_K}$$
et l'on note $\Dlog(V)$ le sous-$\BK$-espace vectoriel de $\D(V\otimes \Q_p[e])$ engendr\'e par $\Dlogp(V)$.

\begin{lem}
\begin{enumerate}[(i)]
\item Le $\BKp$-module $\Dlogp(V)$ est stable sous la d\'erivation $\frac d {de}$. 
\item On a $ \Dlog(V) =  (\BK.\Blogp\otimes_{\Q_p}V)^{H_K}.$
\item Si $V$ est une repr\'esentation de dimension $d$, alors $\Dlog(V)$ est un sous-\phiGamma-module \'etale de $\D(V[d])$ de dimension \'egale au rang de $\Dlogp(V)$ sur $\BKp$ et inf\'erieure ou \'egale \`a $d$. De plus, la repr\'esentation associ\'ee \`a $\Dlog(V)$ est de hauteur finie.
\end{enumerate}
\end{lem}
\Proof
La d\'erivation commute avec l'action de Galois sur $\Blogp\otimes_{\Q_p}V$ :  soit $x = \sum_{k=0}^n a_ke^k$ avec $a_k \in V \otimes \Bp$. Alors
$$ \frac d{de}\tau(x) =  \sum_{k=1}^n \tau(ka_k)(e+1)^{k-1} = \tau\left( \frac d{de}x\right).$$
On obtient ainsi le premier point.
On en d\'eduit imm\'ediatement que $\Dlogp(V)$ est en fait un sous-$\BKp$-module de $\Dp(V[d])$, si bien qu'il est \'egalement de type fini.

En outre, $\Dlogp(V)$ est stable par $\varphi$ ; il suit alors de \cite[Th\'eor\`eme B.1.4.2]{Fontaine_rep_p_corloc} d'une part que $\Dlog(V)$ est \'etale et de hauteur finie, d'autre part que l'application produit
$$ \BK \otimes_{\BKp} \Dlogp(V) \rightarrow \Dlog(V)$$
est injective. On en d\'eduit que la dimension de $\Dlog(V)$ \'egale le $\BKp$-rang de $\Dlogp(V)$ et que $\iota$ induit une injection  $\BK$-lin\'eaire $\Dlog(V) \hookrightarrow \D(V)$, ce qui montre que $\dim_{\BK} \Dlog(V) \le d$.

Montrons enfin que l'injection naturelle 
$$\eta : \Dlog(V) \rightarrow (\BK.\Blogp\otimes_{\Q_p}V)^{H_K}$$
est surjective. Pour cela, on consid\`ere $x = \sum_{i=0}^nb_iw_i \in (\BK.\Blogp\otimes_{\Q_p}V)^{H_K}$ o\`u les $w_i \in \Blogp\otimes_{\Q_p}V$ et les $b_i\in \BK$ sont lin\'eairement ind\'ependants sur $\Bp$. On a alors pour tout $h\in H_K$,
$$h(x) = \sum_{i=0}^nb_ih(w_i) = \sum_{i=0}^nb_iw_i$$
si bien que pour tout $i$, $w_i \in \Dlogp(V)$ et $\eta$ est surjective. Ceci conclut le lemme.
\Eproof

\begin{lem}\label{DptildeisoDlogp}
Si $V$ est une repr\'esentation de dimension $d$, et si $W$ est la sous-repr\'esentation de $V[d]$ correspondant \`a $\Dlog(V)$, alors
\begin{enumerate}
\item  l'injection $\Dp(W) \hookrightarrow \Dlogp(V)$ est un isomorphisme,
\item $V$ est na\"ive si et seulement si $W$ est cristalline de dimension $d$,
\item dans ce cas, $W$ s'identifie \`a la repr\'esentation $\widetilde V$ d\'efinie dans la Proposition \ref{cristallineassociee}.
\end{enumerate}
\end{lem}
\Proof
L'injection $W \hookrightarrow V[d]$ induit une injection 
$$\left(W\otimes \BKp\right)^{H_K} = \Dp(W) \hookrightarrow \Dlogp(V) \simeq \left(V[d]\otimes \BKp\right)^{H_K}.$$
Par ailleurs, $\Dlogp(V)$ est un sous-$\BKp$-module de $\Dlog(V) = \D(W)$ stable par $\varphi$, c'est donc un sous-$\BKp$-module de $\Dp(W)$. Ces deux modules sont donc \'egaux.

Si $W$ est cristalline de rang $d$, la Proposition \ref{cristallineassociee} montre que $V$ est na\"ive.

Supposons $V$ na\"ive.
La sous-repr\'esentation  $\widetilde V$ de $V[d]$ est cristalline de dimension $d$ et s'injecte dans $V\otimes \Q_p[e]$. L'injection $ \Dp(\widetilde V) \hookrightarrow \Dlogp(V)$
montre que le $\BKp$-module $\Dlogp(V)$ est de rang au moins $d= \rg_{\BKp}\Dp(\widetilde V)$ ; et le lemme pr\'ec\'edent montre que ce rang est exactement $d$.
Ainsi, $\D(\widetilde V) = \Dlog(V)$ et $\widetilde V = W$.
\Eproof
\begin{prop}\label{Wachdlog}
Si $V$ est une repr\'esentation $p$-adique de dimension $d$, alors $V$ est na\"ive si et seulement
s'il existe un sous-$\BKp$-module libre $\NN$ de rang $d$ de $\Dlogp(V)$ stable par $\Gamma$ et un entier $h \in \Z$ tels que l'action de $\Gamma$
sur $(\NN/X\NN)(h)$ est triviale.
\end{prop}
\Proof Il suffit de combiner le lemme pr\'ec\'edent avec le Th\'eor\`eme \ref{Wach}. \Eproof

De m\^eme, le Th\'eor\`eme \ref{Berger1} devient
\begin{prop}\label{WachdansDlog}
Si $T$ est un r\'eseau d'une repr\'esentation $V$ na\"ive positive, alors il existe un
unique sous-$\AKp$-module $\Nlog(T)$ de $\Dlogp(T)$ satisfaisant les conditions suivantes :
\begin{enumerate}[(i)]
\item $\Nlog(T)$ est un $\AKp$-module libre de rang $d = \dim_{\Q_p}V$ ;
\item l'action de $\Gamma$ pr\'eserve $\Nlog(T)$ et est triviale sur $\Nlog(T)/X\Nlog(T)$ ;
\item il existe un entier $r\geq 0$ tel que $X^r\Dlogp(T) \subset \Nlog(T)$.
\end{enumerate}
De plus, $\Nlog(T)$ est stable par $\varphi$. Enfin, si l'on pose $\Nlog(V) = \BKp\otimes_{\AKp}\Nlog(T)$, alors
$\Nlog(V)$ est l'unique sous-$\BKp$-module de $\Dlogp(V)$ qui v\'erifie l'analogue des conditions ci-dessus.
\end{prop}

\begin{prop}\label{WachfiltdansDlog}
Si $V$ est une repr\'esentation na\"ive positive, alors $\Nlog(V)$ est stable par l'op\'erateur $N = -\frac 1 X \frac d {de}$.
Si l'on munit $\Nlog(V)$ de la filtration
$$ \Fil i \Nlog(V) = \left\{x \in \Nlog(V), \varphi(x) \in \left(\frac{\varphi(X)}X\right)^i \Nlog(V)\right\},$$
alors l'inclusion $\Dcris(\widetilde V) \subset \BrigKp \otimes_{\BKp}\Nlog(V)$ induit un isomorphisme
$$\Dst(V) \simeq \Nlog(V)/X\Nlog(V)$$ de $\phiN$-modules filtr\'es.
\end{prop}
\Proof Vu le Lemme \ref{DptildeisoDlogp} et le Th\'eor\`eme \ref{Berger2}, on peut travailler sur $\NN(\widetilde V)$
o\`u il suffit de montrer que l'op\'erateur $N$ d\'efini co\"incide modulo $X$ avec la monodromie de $\Dst(V)$.

L'application $\frac d {de}\otimes \varepsilon^{-1}$ induit un morphisme $G_K$-\'equivariant de $V[d]$ vers
$V[d-1](-1) \subset V[d](-1)$, donc de $\widetilde V$ vers $\widetilde V(-1)$, et ainsi de $\NN(\widetilde V)$ vers $\NN(\widetilde V(-1))$
qui n'est autre que $X\NN(\widetilde V)(-1)$. Ainsi, $N = -\frac 1 X \frac d {de}$ agit sur $\NN(\widetilde V)$.
Or, d'apr\`es \cite[Th\'eor\`eme III.4.4]{Berger_limitesrepcris}
$$\widetilde \NN := \NN(\widetilde V) \otimes_{\BKp}\BrigKp = \Dcris(\widetilde V) \oplus X\NN(\widetilde V)$$ d'o\`u
$$-\frac d {de} : \Dcris(\widetilde V) \oplus X\widetilde \NN \rightarrow X\Dcris(\widetilde V) \oplus X^2\widetilde \NN$$
co\"incide avec l'application $tN$ de $\Dcris(\widetilde V)$ dans
$$t\Dcris(\widetilde V) \subset X\Dcris(\widetilde V) \oplus (t-X)\Dcris(\widetilde V)
 \subset X\Dcris(\widetilde V) \oplus X^2\widetilde \NN.$$
On en d\'eduit que l'application $-\frac 1 X \frac d {de} : \NN(\widetilde V) \rightarrow \NN(\widetilde V)$ co\"incide bien
modulo $X$ avec l'op\'erateur $-\frac 1 t \frac d {de}$ de $\Dcris(\widetilde V)$.
\Eproof

\subsection{Module de Wach pour les repr\'esentations na\"ives}
On construit dans ce paragraphe un module de Wach dans $\DLp(V)$ et l'on montre le r\'esultat principal.

Comme $\Dlogp(V) = (\Bp[e]\otimes_{\Q_p}V)^{H_K}$ o\`u l'action de $G_K$ sur $e$ se factorise par $G_L$, on a
$$ \Dlogp(V)= (\DLp(V) \otimes_{\Q_p} \Q_p[e])^{\tau=1}.$$
Nous noterons (certes abusivement) les \'el\'ements de $\Dlogp(V)$ comme des polyn\^omes en $e$ \`a coefficients dans $\DLp(V)$, invariants sous $\tau$.
De plus, $\Dlogp(V)$ est stable par $XN=-\frac d {de}$.

On d\'efinit 
$$\DlogLp(V) = \bigcup_{n\in \N}\ker \left[ (\tau-1)^n : \DLp(V)\rightarrow \DLp(V)\right],$$
le plus grand sous $\BKp$-module de $\DLp(V)$ sur lequel $\tau$ est localement unipotent. La s\'erie $\log\tau$ est alors bien d\'efinie sur $\DlogLp(V)$ et en d\'efinit un endomorphisme localement unipotent. Par ailleurs, il est clair que $\DlogLp(V)$ est stable sous les actions de $G$ et $\varphi$.

\begin{prop}\label{dlogversdl}
Si $V$ est une repr\'esentation $p$-adique de dimension $d$, alors
\begin{enumerate}[(i)]
\item le morphisme de $\BKp$-modules
\begin{eqnarray*}
f: \Dlogp(V) & \rightarrow & \DLp(V)\\
 P(e) & \mapsto & P(0)
\end{eqnarray*}
est injectif , compatible \`a $\Gamma$ et $\varphi$ et d'image $\DlogLp(V)$ ;
\item le morphisme de $\BKp$-modules
\begin{eqnarray*}
g: \DlogLp(V) & \rightarrow & \Dlogp(V)\\
 x & \mapsto & \sum_{n=0}^\infty \frac {e^n} {n!}(-\log \tau)^n x
\end{eqnarray*}
est inverse de $f$, et donc compatible \`a $\Gamma$ et $\varphi$ ;
\item on a $\tau \circ f = f \circ \exp(- \frac d {de})$, et  $f \circ(- \frac d {de}) = (\log \tau) \circ f$ ;
\item si $\Dlogp(V)$ est de rang $d$, alors $\DlogLp(V)$ est un $\BKp$-module libre de rang $d$
qui engendre le $\BL$-espace vectoriel $\DL(V)$.
\end{enumerate}
\end{prop}
\Proof
On note $P^{\sigma}$ le polyn\^ome dont les coefficients sont ceux de $P$ sur lesquels on a agi par $\sigma = \varphi$, $\gamma$ ou $\tau$.
 
Soit $P(e) = \sum_{k=0}^n a_ke^k \in \Dlogp(V)$. Alors pour tout $m \in \Z$, $P^{\tau^m}(e+m) = P(e)$ si bien que
$P(0)=0$ implique $\sum_{k=1}^n a_km^k=0$ pour tout $m \in \Z$, et l'on en d\'eduit ais\'ement que $P =0$.

La compatibilit\'e \`a $\Gamma$ et $\varphi$ provient de ce que
$$ \varphi(P(e)) = P^\varphi(e) \ \textrm{ et } \ \gamma(P(e)) = P^\gamma\left(\frac 1 {\chi(\gamma)} e\right), \ \forall \gamma \in \Gamma.$$
De plus, $ \tau(P(e)) = P^\tau(e+1) = P(e)$ donne $ \tau(P(0)) = P(-1)$. La relation $\tau \circ f = f \circ \exp(- \frac d {de})$ est le
d\'eveloppement de Taylor en $0$ de $P(-1)$. On en d\'eduit que $\tau$ est localement unipotent sur l'image de $f$, donc que cette image est incluse dans $\DlogLp(V)$.

La d\'erivation correspond donc \`a $-\log \tau$. Ainsi pour tout $x \in \Dlogp(V)$,
$$ x = g\circ f(x) = \sum_{n=0}^\infty \frac {e^n} {n!}(-\log \tau)^n f(x)$$
et il est clair que $g$ est inverse de $f$, d\'efini sur $\DlogLp(V)$ qui est donc l'image de $f$.
De plus, si $\Dlogp(V)$ est de rang $d$, alors l'image de $\Dlog(V)$ par l'application induite par l'injection $\iota$ de la Proposition \ref{injectionlogb} est $\D(V)$ car c'en est un sous-$\BK$-espace vectoriel de m\^eme dimension ; donc elle engendre $\DL(V)$. Or
$$ \iota (x) = \sum_{n=0}^\infty \frac {b^n} {n!}(-\log \tau)^n f(x)$$
si bien que l'image de $f$ engendre $\DL(V)$.
\Eproof

\begin{cor} Si $V$ est na\"ive, elle est isomorphe \`a la repr\'esentation cristalline $\widetilde V$ comme repr\'esentation de $G_{\Kpi}$.
\end{cor}
\Proof
Le $(i)$ de la proposition montre que les $\phiG$-modules de $V$ et $\widetilde V$ deviennent isomorphes si l'on oublie l'action de $\tau$.
\Eproof

Remarquons plus g\'en\'eralement que si $\Dlogp(V)$, ou $\DlogLp(V)$, est de rang $d$, $V$ est isomorphe comme repr\'esentation de $G_{\Kpi}$ \`a la repr\'esentation associ\'ee \`a $\Dlog(V)$, qui est de hauteur finie.

On d\'eduit \'egalement de la proposition le r\'esultat suivant.
\begin{cor}\label{produitinjectif}
L'application produit $\BL\otimes_{\BKp}\DlogLp(V) \rightarrow \DL(V)$ est injective. Son image est un sous-\phiG-module \'etale sur $\BL$ de $\DL(V)$.
\end{cor}
\Proof
Il s'agit premi\`erement de v\'erifier que si $(x_i)_i$ est une base de $\DlogLp(V)$, cette famille reste libre sur $\BL$ dans $\DL(V)$. 
Commen\c cons par montrer qu'elle est libre sur $\BK$. Si l'on a une relation de liaison de $(x_i)_i$ sur $\BK$, on peut lui appliquer une puissance convenable de $\tau-1$ pour obtenir une relation de liaison sur $\BK$ d'\'el\'ements de $\Dp(V)$ libres sur $\BKp$. Ceci n'est pas possible du fait de \cite[Th\'eor\`eme B.1.4.2]{Fontaine_rep_p_corloc}.

L'application $\tilde h : x \mapsto \sum_{n=0}^\infty \frac {b^n} {n!}\otimes(-\log \tau)^nx$ d\'efinit un changement de base $x_i \mapsto \tilde h(x_i)$ de $\BL \otimes \DlogLp(V)$.
Or  $h : x \mapsto \sum_{n=0}^\infty \frac {b^n} {n!}(-\log \tau)^nx$ d\'efinit d'apr\`es les Proposition \ref{dlogversdl} et Proposition \ref{injectionlogb} une injection $\DlogLp(V) \hookrightarrow \D(V)$. 
On est donc ramen\'e \`a v\'erifier que la famille $(h(x_i))_i$ est libre sur $\BL$. Ceci est une cons\'equence imm\'ediate de l'isomorphisme de comparaison
$$ \B \otimes \D(V) \rightarrow \B \otimes V$$
qui montre que le produit $ \BL \otimes \D(V) \rightarrow \DL(V)$ est un isomorphisme.

Enfin,  \cite[Th\'eor\`eme B.1.4.2]{Fontaine_rep_p_corloc} montre \'egalement que $\BL \Dp(V)$ est \'etale et en prenant un r\'eseau $T$ de $V$, on montre par une r\'ecurrence imm\'ediate que pour tout $n$, 
$\varphi\left(\BL(\tau-1)^{-n}(\Dp(T))\right)$ engendre $\BL(\tau-1)^{-n}(\Dp(T))$ sur $\BL$, si bien que $\BL\DlogLp(V)$ est bien \'etale.
\Eproof

\begin{lem}
Si $M$ est un sous $\AKp$-module de $\AL^d$ de type fini, stable par $\tau$, alors pour tout $k$ assez grand, $(\tau-1)^kM \subset pM$.
\end{lem}
\Proof
Comme $M$ est un $\AKp$-module de type fini, c'est l'image d'un morphisme $(\AKp)^s \rightarrow \AL^d$ continu pour la topologie $p$-adique. On v\'erifie alors ais\'ement que pour $m$ suffisamment grand $p^m\AL^d \cap M \subset pM$.

Fixons $(e_i)_{1\le i \le d}$ la base canonique de $\AL^d$ et $(m_j)_{1 \le j \le r}$ une famille g\'en\'eratrice de $M$ sur $\AKp$ que l'on \'ecrit dans la base canonique
$$\forall 1\le j \le r, \ m_j=\sum_{i=1}^d x_{ij}e_i.$$
On est ramen\'e \`a trouver un $k$ tel que pour tout $1\le i \le d$, $ 1\le j \le r$, $(\tau-1)^kx_{ij} \in p^m\AL$. Mais tout $x\in \AL$ s'\'ecrit $x=\sum_{l\in \N}p^l[x_l]$ avec $x_l \in \EL$ o\`u $\EL = \E^{G_L} = \cup_{n\in \N}\EL^{(\tau^{p^n}=1)}$ d'apr\`es la th\'eorie du corps des normes \cite[Th\'eor\`eme 3.2.2.]{WinCorNorm}. 
Il est alors clair que $k=p^r$ avec $r$ assez grand convient.
\Eproof

\begin{lem}\label{ssmoddeDlogLp}
Si $V$ est une repr\'esentation et si $M$ est un sous-$\BKp$-module de type fini de $\DL(V)$ stable par $G$, alors $\tau$ est unipotent sur $M$. 

En particulier, si $M$ est un sous-$\BKp$-module de type fini de $\DLp(V)$ stable par $G$, alors $M\subset \DlogLp(V)$. 
\end{lem}
\Proof
Fixons $V_0$ un r\'eseau de $V$ et $M_0 = M\cap \DL(V_0)$. Comme $\DL(V_0)$ poss\`ede une base fixe par $\tau$, celle de $D(V_0)$, on peut l'identifier (si on le consid\`ere muni de la seule action de $\tau$) \`a $\AL^d$ et appliquer le lemme pr\'ec\'edent : il existe $k$ tel que $(\tau-1)^kM_0 \subset pM_0$.

On en d\'eduit que la s\'erie $\log \tau$ converge vers un endomorphisme de $M$. Soit $\Delta \in M_d(\BKp)$ sa matrice dans une base $\mathcal B$ de $M_0$ et $T \in GL_d(\AKp)$ celle de $\tau$.
Pour tout $\gamma \in \Gamma$, notons $H(\gamma) \in GL_d(\AKp)$ sa matrice dans $\mathcal B$. Alors la relation de commutation $H(\gamma)T^{\gamma} = T^{\chi(\gamma)}H(\gamma)$ donne
$$H(\gamma)\Delta^{\gamma}H(\gamma)^{-1} = \chi(\gamma)\Delta.$$
Si l'on \'ecrit $\chi_\Delta$ le polyn\^ome caract\'eristique de $\Delta$ sous la forme $\chi_\Delta(\lambda) = \lambda^d+a_{d-1}\lambda^{d-1} + \dots + a_0$, alors pour tout $0\le i \le d-1$ et pour tout $\gamma \in \Gamma$, on a $\gamma(a_i) = \chi(\gamma)^{d-i}a_i$. On en d\'eduit que tous les $a_i$ sont nuls, puisque $a_i \in \BKp \subset (\Bcris)^{H_K}$ et qu'il est bien connu que les solutions de $\gamma(x) = \chi(\gamma)^{k}x$ dans $(\Bcris)^{H_K}$ sont les $x=x_0t^k$ avec $x_0 \in K$, d'o\`u pour tout $i$, $a_i \in Kt^{d-i}\cap \BKp = \{0\}$.

Il suit que $\Delta$ est nilpotent, et ainsi $\tau$ est unipotent.
\Eproof

\begin{prop}\label{caracDlogLp}
Si $V$ est une repr\'esentation $p$-adique de dimension $d$, alors $\DlogLp(V)$ est le plus grand sous-$\BKp$-module $M$ de type fini de $\DLp(V)$ stable par $G$. Si $V$ est na\"ive, $\DlogLp(V)$ est libre de rang $d$ et engendre $\DL(V)$ sur $\BL$.
\end{prop}
\Proof Le module $\Dlogp(V)$ est de rang $d$ d'apr\`es le Lemme \ref{DptildeisoDlogp}, donc $\DlogLp(V)$ aussi
et la proposition se d\'eduit de la Proposition \ref{dlogversdl} et du Lemme \ref{ssmoddeDlogLp}.
\Eproof

\begin{thm}\label{Wachinitsst}
Si $V$ est une repr\'esentation $p$-adique de dimension $d$, alors $V$ est na\"ive si et seulement
s'il existe un sous-$\BKp$-module $\NN$ libre de rang $d$ de $\DLp(V)$, stable par $G$,
et un entier $h \in \Z$ tels que l'action de $\Gamma$ sur $(\NN/X\NN)(h)$ est triviale.
\end{thm}
\Proof C'est une cons\'equence de la Proposition \ref{Wachdlog}, de la Proposition \ref{dlogversdl} et du Lemme \ref{ssmoddeDlogLp}.
\Eproof

\begin{thm}\label{Wachsst}
Si $V$ est une repr\'esentation na\"ive positive de dimension $d$, alors il existe un unique
sous-$\BKp$-module $\NN(V)$ de $\DLp(V)$ satisfaisant les conditions suivantes :
\begin{enumerate}
\item $\NN(V)$ est un $\BKp$-module libre de rang $d$ qui engendre $\DL(V)$;
\item l'action de $G$ pr\'eserve $\NN(V)$, et est triviale sur $\NN(V)/X\NN(V)$ ;
\item il existe un entier $r\geq 0$ tel que $X^r\DlogLp(V) \subset \NN(V)$.
\end{enumerate}
De plus, $\NN(V)$ est stable par $\varphi$.
\end{thm}
\Proof
Ce th\'eor\`eme, except\'e le fait que $(\tau-1)(\NN(V)) \subset X\NN(V)$ est une cons\'equence de la Proposition \ref{WachdansDlog}, de la Proposition \ref{dlogversdl} et du Lemme \ref{ssmoddeDlogLp}.

La Proposition \ref{WachfiltdansDlog} montre que $\Nlog(V)$ est stable par $\frac 1 X \frac d{de}$ donc que $\frac d {de} (\Nlog(V)) \subset X \Nlog(V)$, on en d\'eduit que $\log(\tau)(\NN(V)) \subset X\NN(V)$ si bien que $\tau$ est trivial sur $\NN(V)/X\NN(V)$.
\Eproof

\begin{thm}\label{Wachsstfiltre}
Si $V$ est une repr\'esentation na\"ive positive, la somme $\frac 1 X \log \tau$ converge vers un op\'erateur $N$ nilpotent sur $\NN(V)$.
Si l'on munit de plus  $\NN(V)$ de la filtration
$$ \Fil i \NN(V) = \left\{x \in \NN(V), \varphi(x) \in \left(\frac{\varphi(X)}X\right)^i \NN(V)\right\},$$
alors il existe un isomorphisme
$\Dst(V) \simeq \NN(V)/X\NN(V)$ de $\phiN$-modules filtr\'es.
\end{thm}
\Proof
C'est une cons\'equence des Proposition \ref{WachfiltdansDlog} et \ref{dlogversdl} et du Lemme \ref{ssmoddeDlogLp}.
\Eproof

\begin{rmq} \label{rmqphimodfiltre}
\begin{enumerate}[(i)]
\item  Il semble que ce soit un
ph\'enom\`ene g\'en\'eral pour les repr\'e\-sen\-ta\-tions semi-stables que la monodromie soit encod\'ee par la partie en $\log \tau$ de l'alg\`ebre de Lie de $G$.
\item Comme il a \'et\'e dit, $\tau$ ne d\'epend que tu choix de $\varepsilon$ ; si l'on change $\varepsilon$ en $\varepsilon^\alpha$,
$\tau$ est chang\'e en $\tau^\alpha$ et l'op\'erateur $N = \frac 1 X \log \tau$ n'est pas chang\'e modulo $X$.
De plus, $V$ est cristalline si et seulement si $N=0$, c'est-\`a-dire $\tau=1$, et alors $\NN(V)\subset D(V)$ est le module de Wach du Th\'eor\`eme
\ref{Berger1}.
\item Les modules $\NN(\widetilde V)$, $\Nlog(V)$ et $\NN(V)$  sont isomorphes comme $\varphi$-modules filtr\'es.
\end{enumerate}
\end{rmq}

\subsection{Modules de Wach dans $\D(V)$}\label{sectionmoche}

On peut encore, \`a l'aide des r\'esultats pr\'e\-c\'e\-dents, construire un module de Wach pour les repr\'esentations na\"ives
dans le $\phiGamma$-module $\D(V)$. Toutefois il d\'epend d'un choix, celui de l'\'el\'ement $b$ du Lemme \ref{elementb}.
On fixe donc $b \in \BL$ tel que $(h-1)b = \xi(h)$ pour tout $h \in H_K$, ou de mani\`ere \'equivalente
$(\tau-1)b=1$. On remarque que $(\varphi - 1)b, (\chi(\gamma)\gamma - 1)b \in \BK$ pour $\gamma \in \Gamma<G$.
\footnote{Si $p \neq 2$, $\Gamma$ est procyclique de g\'en\'erateur $\gamma$
et les couples $((\varphi - 1)b\otimes \varepsilon, (\gamma - 1)(b\otimes \varepsilon))$
correspondent \`a l'image de $p$ par l'application de Kummer sur le premier groupe de cohomologie du complexe de Herr qui calcule la cohomologie galoisienne de $\Q_p(1)$. Pour nos constructions, le choix de $b$ est \'equivalent au choix d'un tel couple.}

On d\'efinit ensuite $\BLlog$ (respectivement $\BLlogp$) comme la plus petite sous $\BLp$-alg\`ebre de $\BL$ contenant $\BLp[b]$ (respectivement $b$),
stable par $\varphi$ et $G$.
On d\'efinit de m\^eme $\BKlog$ (respectivement $\BKlogp $) comme la plus petite sous $\BKp$-alg\`ebre de $\BK$ contenant $\BKp[(\varphi - 1)b, (\chi(\gamma)\gamma - 1)b, \gamma \in \Gamma]$
(respectivement $(\varphi - 1)b, (\chi(\gamma)\gamma - 1)b, \gamma \in \Gamma$),
stable par $\varphi$ et $\Gamma$.
On a alors
$$ \BLlog = \BLp \oplus \BLlogp \ \textrm{ et } \BKlog = \BKp \oplus \BKlogp.$$

\begin{lem}
On a :
$$ \BKlog = \left(\BLlog\right)^{\tau=1} \ \textrm{ et }   \BKlogp = \left(\BLlogp\right)^{\tau=1}.$$
\end{lem}
\Proof
Les deux \'egalit\'es se traitent de m\^eme mani\`ere, consid\'erons la premi\`ere. Vu le Lemme \ref{elementb}, il suffit de montrer que $\BLlog = \BKlog.\BLp[b]$.
Or, il est clair que $ \BKlog.\BLp[b]\subset\BLlog $ et que $ \BKlog.\BLp[b]$ est stable par $\varphi$ et $\Gamma$, d'o\`u l'\'egalit\'e.
\Eproof

\begin{lem}
Si $V$ est une repr\'esentation na\"ive positive et si $M $ est le $\BKlog$-module $ \left(\BLlog\otimes \DLp(V)\right)^{\tau=1}$ engendr\'e par l'image de l'application $\iota$ de la
Proposition \ref{injectionlogb}, alors $\iota$ induit un isomorphisme de $\BKp$-modules compatible \`a $\varphi$ et $\Gamma$
$$ \Dlogp(V) \simeq M / (\BKlogp M).$$
\end{lem}
\Proof
Cela provient de ce que $\BKp \simeq \BKlog/\BKlogp$.
\Eproof

On en d\'eduit les propositions suivantes.

\begin{prop}
Si $V$ est une repr\'esentation na\"ive positive de dimension $d$, soit $M = \left(\BLlog\otimes \DLp(V)\right)^{\tau=1}$, alors il existe un unique
sous-$\BKlog$-module $\NN_K(V)$ de $\D(V)$ satisfaisant les conditions suivantes :
\begin{enumerate}
\item $\NN_K(V)$ est un $\BKlog$-module libre de rang $d$ ;
\item l'action de $\Gamma$ pr\'eserve le module $\NN_K(V)$ et est triviale sur le quotient $\NN_K(V)/(X,\BKlogp)\NN_K(V)$ ;
\item il existe un entier $r\geq 0$ tel que $X^rM \subset \NN_K(V)$.
\end{enumerate}
De plus, $\NN_K(V)/\BKlogp\NN_K(V)$ est stable par $\varphi$.
\end{prop}

\begin{prop}
Si $V$ est une repr\'esentation semi-stable positive, et si l'on munit $M=\NN_K(V)/\BKlogp\NN_K(V)$ de la filtration
$$ \Fil i M = \left\{x \in M, \varphi(x) \in \left(\frac{\varphi(X)}X\right)^i M\right\},$$
 alors il existe un isomorphisme
$\Dst(V) \simeq M/XM$ de $\varphi$-modules filtr\'es.
\end{prop}

\section{Applications}

Les foncteurs d'inclusion de la cat\'egorie des repr\'esentations respectivement de hauteur finie, cristallines et na\"ives de $G_K$ vers $\Rep_{G_K}$ admettent des adjoints \`a droite $V \mapsto V^{\hf}$, $V \mapsto V^{\cris}$ et $V \mapsto V^{\naif}$ qui \`a une repr\'esentation $V$ associent ses plus grandes sous-repr\'esentations respectivement de hauteur finie, cristalline et na\"ive. On commence par la remarque suivante.
\begin{lem}\label{hfcris}
Si $V$ est une repr\'esentation semi-stable de $G_K$, alors $V^{\cris}=V^{\hf}$ et  $\Dp(V) \simeq \Dp(V^{\cris})$ comme \phiGamma-modules sur $\BKp$.
\end{lem}
\Proof
Le  Th\'eor\`eme B.1.4.2 de \cite{Fontaine_rep_p_corloc} montre que pour toute repr\'e\-sen\-ta\-tion $V$ de $G_K$, $\BK\Dp(V)\subset D(V)$ est un sous-\phiGamma module \'etale de $D(V)$, il correspond donc \`a la plus grande sous-repr\'esentation de hauteur finie de $V$. Si $V$ est semi-stable, le Th\'eor\`eme de \cite{Wach_pot_cris} montre que cette repr\'esentation est en fait potentiellement cristalline, donc cristalline. \Eproof

\begin{prop}\label{propappli}
Si $V$ est une repr\'esentation de $G_K$ telle que $V^{\cris}$ est positive, $M \subset \Dp(V^{\cris}) \subset \Dp(V)$ est le module de Wach de $V^{\cris}$ et $n$ est l'indice de nilpotence de $\log \tau$ sur $\DlogLp(V)$ alors le $\BKp$-module $\NN$ engendr\'e par  $\{X^{i}(\log \tau)^{-i}(M) \ ; \ 0 \le i \le n\}$ satisfait les propr\i\'et\'es :
\begin{enumerate}
\item $\NN$ est un $\BKp$-module libre et $\BL\otimes_{\BKp}\NN  \rightarrow \DL(V)$ est injectif ;
\item l'action de $G$ pr\'eserve $\NN$, et elle est triviale sur $\NN/X\NN$ ;
\item il existe un entier $r\ge n-1$ tel que $X^r(\log \tau)^{-n}(\Dp(V^{\cris})) \subset \NN$.
\end{enumerate}
Enfin, $\BL\NN$ est le sous-\phiG-module \'etale $\DL(V^{\naif})$ de $\DL(V)$ et $\NN$ est son unique sous-$\BKp$-module satisfaisant les propri\'et\'es du Th\'eor\`eme \ref{Wachsst}.
\end{prop}
\Proof 
Le $\BKp$-module $\NN$ est libre comme sous-module de $\DlogLp(V)$ et l'injectivit\'e de $\BL \otimes_{\BKp} \NN  \rightarrow \DL(V)$ suit du Corollaire \ref{produitinjectif}. On sait que $G$ pr\'eserve $M$ et est trivial sur $M/XM$. Si $(\log \tau)^i(x) \in M$, et $g\in G$, alors
$$ g \left((\log \tau)^i(x)\right) = (\chi(g))^i(\log \tau)^i(g(x)) =y$$
o\`u $y\in M$ et $y-(\log\tau)^i(x) \in XM$.  En particulier, $g(x) \in (\log\tau)^{-i}M$ et $X^ig(x) \in \NN$. On d\'eduit alors le second point de :
$$ (\log \tau)^i(g-1)(X^ix)  = X^i \left(\left( \frac{g(X)}X\right)^i \frac y{\chi(g)^i} - (\log \tau)^i(x)\right) \in X^{i+1}M.$$
Le troisi\`eme point suit de l'existence d'un entier $r'$ tel que $X^{r'}\Dp(V^{\cris}) \subset M$.

 La preuve que $\BL\NN$  est \'etale est la m\^eme que dans le Corollaire \ref{produitinjectif}.
 
On en d\'eduit que $V$ poss\`ede une sous-repr\'esentation $W$ telle que $\DL(W) = \BL\NN = \BL(\log \tau)^{-n}(\Dp(V^{\cris}))$. 
La seconde propri\'et\'e, combin\'ee avec le Th\'eo\-r\`e\-me \ref{Wachinitsst}, montre que $W$ est na\"ive. 
Si maintenant $W' \subset V$ est na\"ive, soit $\NN' \subset \DlogLp(V)$ son module de Wach au sens du Th\'eor\`eme \ref{Wachsst}, alors $\NN' \subset \cup_{k\in \N}(\tau-1)^{-k}\Dp(V^{\cris}) = (\log \tau)^{-n}(\Dp(V^{\cris}))$ si bien que $W' \subset W$ et on a donc bien $W=V^{\naif}$.
 \Eproof

\begin{cor}
Si $V^{\cris}=V^{\hf}$, alors $\DlogLp(V) = \DlogLp(V^{\naif})$, autrement dit,  $V^{\naif} = (\B\otimes_{\BKp}\DlogLp(V))^{\varphi=1}$.
\end{cor}
\Proof
Si $V^{\cris}=V^{\hf}$, alors $\Dp(V) = \Dp(V^{\cris})$ si bien que $\DlogLp(V) = \cup_{k\in \N}(\tau-1)^{-k}\Dp(V^{\cris})$ et que l'\'enonc\'e est une cons\'equence imm\'ediate de la Proposition \ref{propappli}.
\Eproof

Appelons repr\'esentation (respectivement semi-stable) de $G$ toute repr\'e\-sen\-ta\-tion (respectivement semi-stable) de $G_K$ dont l'action de Galois se factorise par $G$.

\begin{cor}\label{RepGnaives}
Les repr\'esentations semi-stables de $G$ sont na\"ives. Une telle re\-pr\'e\-sen\-ta\-tion est cristalline \ssi elle est fix\'ee par $\tau$.
\end{cor}
\Proof
Si $V$ est fixe par $G_L$, alors $\DL(V) = \BL\otimes V$ et d'apr\`es le Lemme \ref{ssmoddeDlogLp}, $\BKp\otimes V \subset \DlogLp(V) \subset \BLp\otimes V$. Il suit alors du Corollaire \ref{produitinjectif} que $\BKp\otimes V$ est un sous-$\BKp$-module de $\DlogLp(V)$ de m\^eme rang, et comme $\BKp$ est principal, il existe une base $(d_i)_i$ de  $\BKp\otimes V$ et des $\lambda_i \in \BKp$ tels que $(\frac 1 {\lambda_i}d_i)_i$ est une base de  $\DlogLp(V)$. Comme $\Frac \BKp \cap \BLp = \BKp$, les $\lambda_i$ sont inversibles, bref $ \DlogLp(V) = \BKp\otimes V$.

D'apr\`es le Lemme \ref{hfcris}, si $V$ est semi-stable, alors $V^{\cris}=V^{\hf}$, et donc d'apr\`es le corollaire pr\'ec\'edent, $V$ est na\"ive.

Enfin, $V$ est cristalline \ssi $\log \tau=0$ sur $\NN(V)$, ceci \'equivaut ici \`a ce que $\tau=1$ sur $\BKp\otimes V$ donc sur $V$.
\Eproof

Le th\'eor\`eme suivant a \'et\'e conjectur\'e par Breuil \cite{Breuil_CorpsNormes} et prouv\'e par Kisin \cite{Kisin_crys}. On en obtient ici une preuve, dans le cas o\`u $K$ est non ramifi\'e, comme cons\'equence du corollaire pr\'ec\'edent. Notons que les arguments de cette preuve peuvent \^etre r\'eduits et g\'en\'eralis\'es de sorte \`a obtenir une preuve du th\'eor\`eme de Kisin y compris dans le cas ramifi\'e \cite{BeilinsonTavares}.

\begin{thm}
Le foncteur restriction de la cat\'egorie des repr\'esentations cristallines de $G_K$ vers la cat\'egorie des repr\'esentations de $G_{\Kpi}$ est pleinement fid\`ele.
\end{thm}
\Proof
Soient $V_1$ et $V_2$ deux repr\'esentations cristallines, alors la repr\'e\-sen\-ta\-tion $V:=\Hom(V_1, V_2)$ est cristalline. Un \'el\'ement $ v\in V^{G_{\Kpi}}$ appartient en particulier \`a la sous-repr\'esentation $V^{G_L}$ qui est cristalline, donc fix\'ee par $\tau$, ainsi $v\in V^{G_K}$.
\Eproof

Nous allons maintenant d\'ecrire les repr\'esentations semi-stables de $G$. Donnons tout d'abord un exemple.

\begin{exe}
Rappelons que l'on note $V_i=\Sym^{i-1} V_2$ o\`u $V_2 = \Q_p \oplus \Q_p e$ est la repr\'e\-sen\-ta\-tion associ\'ee \`a la courbe de Tate. On a alors l'identification $V_i \simeq \Q_p[e]_{\le i-1}$ avec l'espace des polyn\^omes de degr\'e strictement plus petit que $i$. On note de m\^eme $V_1=\Q_p$. Si $j\ge 0$, le module de Wach de $V_i(-j)$ est le $\BKp$-module libre de base
$ X^{k+j}\otimes e^k$, $0 \le k \le i-1$.
\end{exe}

Le r\'esultat suivant est d\^u \`a Perrin-Riou \cite[Lemme 3.3.4.]{PerrinRiou_IwasawaCorpsLocal}  (cf. aussi \cite[Lemme 3.3.8.]{BergerBpaires}).

\begin{lem}
Si $V$ est une repr\'e\-sen\-ta\-tion de $\Gamma$, alors $V$ est cristalline \ssi $V=\oplus_{j\in \Z}V^{\Gamma=\chi^j}$. 
\end{lem}

\begin{prop}\label{descriptionsstG}
Si $V$ est une repr\'esentation semi-stable de $G$, alors il existe des entiers $n_{i,j}$ presque tous nuls tels que $V\simeq \oplus_{i \ge 1, j\in \Z} V_i(j)^{n_{i,j}}$.
\end{prop}
\Proof
D'apr\`es le lemme pr\'ec\'edent, il existe des entiers $m_j\ge 0$ presque tous nuls tels que :
$$V^{\cris} = V^{\tau=1} \simeq \oplus_{j\in \Z} \Q_p(j)^{m_j}$$
de plus, si l'on tord par une puissance assez grande de $\chi^{-1}$, on peut supposer $V^{\cris}$ positive, c'est-\`a-dire $m_j=0$ pour tout $j > 0$. D'apr\`es la Corollaire \ref{RepGnaives} $V$ est na\"ive, on peut donc  appliquer la Proposition \ref{propappli} qui donne le r\'esultat.
\Eproof

\subsection*{Remerciements}
Une premi\`ere version de ce palimpseste a \'et\'e \'ecrite \`a l'UMPA, laboratoire de math\'ematiques de l'\'ENS de Lyon. J'en remercie les membres de l'\'equipe de th\'eorie des nombres et particuli\`erement Laurent Berger pour des conversations profitables.

\def\cprime{$'$}
\providecommand{\bysame}{\leavevmode\hbox to3em{\hrulefill}\thinspace}
\providecommand{\MR}{\relax\ifhmode\unskip\space\fi MR }
\providecommand{\MRhref}[2]{%
  \href{http://www.ams.org/mathscinet-getitem?mr=#1}{#2}
}
\providecommand{\href}[2]{#2}


\begin{thebibliography}{BTR12}

\bibitem[Ber02]{Berger_equadiff}
Laurent Berger, \emph{Repr\'esentations {$p$}-adiques et \'equations
  diff\'erentielles}, Invent. Math. \textbf{148} (2002), no.~2, 219--284.
  \MR{MR1906150 (2004a:14022)}

\bibitem[Ber04]{Berger_limitesrepcris}
\bysame, \emph{Limites de repr\'esentations cristallines}, Compos. Math.
  \textbf{140} (2004), no.~6, 1473--1498. \MR{MR2098398 (2006c:11138)}

\bibitem[Ber08a]{BergerBpaires}
\bysame, \emph{Construction de {$(\phi,\Gamma)$}-modules: repr\'esentations
  {$p$}-adiques et {$B$}-paires}, Algebra Number Theory \textbf{2} (2008),
  no.~1, 91--120. \MR{2377364 (2009j:14025)}

\bibitem[Ber08b]{Berger_EquadiffModulesfiltres}
\bysame, \emph{\'{E}quations diff\'erentielles {$p$}-adiques et
  {$(\phi,N)$}-modules filtr\'es}, Ast\'erisque (2008), no.~319, 13--38,
  Repr{\'e}sentations $p$-adiques de groupes $p$-adiques. I.
  Repr{\'e}sentations galoisiennes et $(\phi,\Gamma)$-modules. \MR{MR2493215
  (2010d:11056)}

\bibitem[Bre97]{Breuil_Griffiths}
Christophe Breuil, \emph{Repr\'esentations {$p$}-adiques semi-stables et
  transversalit\'e de {G}riffiths}, Math. Ann. \textbf{307} (1997), no.~2,
  191--224. \MR{MR1428871 (98b:14016)}

\bibitem[Bre99]{Breuil_CorpsNormes}
\bysame, \emph{Une application de corps des normes}, Compositio Math.
  \textbf{117} (1999), no.~2, 189--203. \MR{MR1695849 (2000f:11157)}

\bibitem[BTR12]{BeilinsonTavares}
Alexander Beilinson and Floric Tavares~Ribeiro, \emph{On a theorem of {K}isin},
  arXiv:1205.2913v2 (2012).

\bibitem[CF00]{ColmezFontaine_ConstructionSst}
Pierre Colmez and Jean-Marc Fontaine, \emph{Construction des repr\'esentations
  {$p$}-adiques semi-stables}, Invent. Math. \textbf{140} (2000), no.~1, 1--43.
  \MR{MR1779803 (2001g:11184)}

\bibitem[Col99]{ColmezCristallinesHauteurFinie}
Pierre Colmez, \emph{Repr\'esentations cristallines et repr\'esentations de
  hauteur finie}, J. Reine Angew. Math. \textbf{514} (1999), 119--143.
  \MR{MR1711279 (2001h:11147)}

\bibitem[Col02]{Colmez_Banachfinie}
\bysame, \emph{Espaces de {B}anach de dimension finie}, J. Inst. Math. Jussieu
  \textbf{1} (2002), no.~3, 331--439. \MR{MR1956055 (2004b:11160)}

\bibitem[Fon90]{Fontaine_rep_p_corloc}
Jean-Marc Fontaine, \emph{Repr\'esentations {$p$}-adiques des corps locaux.
  {I}}, The Grothendieck Festschrift, Vol.\ II, Progr. Math., vol.~87,
  Birkh\"auser Boston, Boston, MA, 1990, pp.~249--309. \MR{MR1106901
  (92i:11125)}

\bibitem[Fon94]{Fontaine_cor_perio}
\bysame, \emph{Le corps des p\'eriodes {$p$}-adiques}, Ast\'erisque (1994),
  no.~223, 59--111, With an appendix by Pierre Colmez, P\'eriodes $p$-adiques
  (Bures-sur-Yvette, 1988). \MR{MR1293971 (95k:11086)}

\bibitem[Fon04]{Fontaine_ArithmetiqueRepGalois}
\bysame, \emph{Arithm\'etique des repr\'esentations galoisiennes
  {$p$}-adiques}, Ast\'erisque (2004), no.~295, xi, 1--115, Cohomologies
  $p$-adiques et applications arithm{\'e}tiques. III. \MR{MR2104360
  (2005i:11074)}

\bibitem[Kis06]{Kisin_crys}
Mark Kisin, \emph{Crystalline representations and {$F$}-crystals}, Algebraic
  geometry and number theory, Progr. Math., vol. 253, Birkh\"auser Boston,
  Boston, MA, 2006, pp.~459--496. \MR{MR2263197 (2007j:11163)}

\bibitem[Mok98]{Mokrane_Ordinarite}
A.~Mokrane, \emph{Quelques remarques sur l'ordinarit\'e}, J. Number Theory
  \textbf{73} (1998), no.~2, 162--181. \MR{MR1658011 (2000b:11131)}

\bibitem[PR94]{PerrinRiou_IwasawaCorpsLocal}
Bernadette Perrin-Riou, \emph{Th\'eorie d'{I}wasawa des repr\'esentations
  {$p$}-adiques sur un corps local}, Invent. Math. \textbf{115} (1994), no.~1,
  81--161, With an appendix by Jean-Marc Fontaine. \MR{1248080 (95c:11082)}

\bibitem[TR11]{Tavares_AIF}
Floric Tavares~Ribeiro, \emph{An explicit formula for the {H}ilbert symbol of a
  formal group}, Ann. Inst. Fourier (Grenoble) \textbf{61} (2011), no.~1,
  261--318. \MR{2828131}

\bibitem[Wac96]{Wach_pot_cris}
Nathalie Wach, \emph{Repr\'esentations {$p$}-adiques potentiellement
  cristallines}, Bull. Soc. Math. France \textbf{124} (1996), no.~3, 375--400.
  \MR{MR1415732 (98b:11119)}

\bibitem[Win83]{WinCorNorm}
Jean-Pierre Wintenberger, \emph{Le corps des normes de certaines extensions
  infinies de corps locaux; applications}, Ann. Sci. \'Ecole Norm. Sup. (4)
  \textbf{16} (1983), no.~1, 59--89. \MR{MR719763 (85e:11098)}

\end{thebibliography}
\end{document}